\documentclass[twoside,12pt]{article}
\usepackage{CJK,indentfirst,amsmath,amsfonts,amssymb,amsthm,cite}
\usepackage{fancyhdr}
\usepackage[dvips]{graphics}
\usepackage[dvips]{color}
\begin{document}
\vspace*{0.2cm}
\begin{center}\baselineskip=11mm
{\LARGE\bf {M-estimation in high-dimensional linear model}}\\[2em]
\normalsize\bf
Kai Wang$^\dagger$, Yanling Zhu$^{\dagger*}$ \\[0.3cm]

%-----------------footnotes-------------------------------------------------------------------
\begin{figure}[b]
\footnotesize
\rule[-2.5truemm]{5cm}{0.1truemm}\\[1mm]\baselineskip=5mm
{$^\dagger$ School of Statistics and Applied Mathematics, Anhui University of Finance and Economics,\\ Bengbu
233030, \ P.R. CHINA \\
$^*$Corresponding author.\\
 E-mail addresses\,: zhuyanling99\textit{\char64}126.com,\; wangkai050318\textit{\char64}163.com}\\
Supported by the NSF of Anhui Province (No.1508085QA13,1708085MA17).\\

\end{figure}
\end{center}
%-------------------Abstract------------------------------------------------------------------
 \noindent\hrulefill \newline
\begin{center}
\begin{minipage}{15cm}\baselineskip=10mm
\noindent{\bf Abstract}\ \ We mainly study the M-estimation method for the high-dimensional linear regression model, and discuss the properties of M-estimator when the penalty term is the local linear approximation. In fact, M-estimation method is a framework, which covers the methods of the least absolute deviation, the quantile regression, least squares regression and Huber regression. We show that the proposed estimator possesses the good properties by applying certain assumptions. In the part of numerical simulation, we select the appropriate algorithm to show the good robustness of this method.

%-------------------Keywords------------------------------------------------------------------
\vspace{2mm}\noindent {{\bf Keywords:}}{ M-estimation; High-dimensionality; Variable selection; Oracle property; Penalized method}\\[2mm]
{\it 2000 AMS Classification}: 62F12, 62E15, 62J05.\\
\end{minipage}
\end{center}

%\vspace{0.1cm}
\noindent\hrulefill \newline
%---------------------------Text-------------------------------------------------------------
\baselineskip=10mm
\section*{\normalsize\bf 1 \quad{\normalsize\bf Introduction} }
For the classical linear regression model $Y=X\beta+\varepsilon$, we are interested in the problem of variable selection and estimation, where $Y=(y_1,y_2,...,y_n)^T$ is the response vector, $X=(X_1,X_2,...,X_{p_n})=(x_1,x_2,...,x_n)^T=(x_{ij})_{n\times p_n}$ is an $n\times p_n$ design matrix,
 and $\varepsilon=(\varepsilon_1,\varepsilon_2,...,\varepsilon_n)^T$ is a random vector. The main topic is how to estimate the coefficients vector $\beta\in \textrm{R}^p_n$ when $p_n$ increases with sample size $n$ and many elements of $\beta$ equal zero. We can transfer this problem into a minimization of a penalized least squares objective function
 $$\hat{\beta_n}=\text{arg}\min_\beta Q_n(\beta_n),\ Q_n(\beta_n)=\|Y-X\beta_n\|^2+\sum_{j=1}^{p_n}p_{\lambda_n}(|\beta_nj|),$$
 where $||\cdot||$ is the $l_2$ norm of the vector,\ $\lambda_n$ is a tuning parameter, and $p_{\lambda_n}(|t|)$ a penalty term.
We have known that least squares estimation is not robust, especially when the data exists abnormal values or the error term has the heavy tailed distribution.
\par In this paper we consider the loss function be least absolute deviation,i.e., minimize the following objective function:
\[\hat{\beta_n}=\text{arg}\min_\beta Q_n(\beta_n),\ Q_n(\beta_n)=\frac{1}{n}\sum_{i=1}^{n}|y_i-x_i^T\beta_n|+\sum_{j=1}^{p_n}p_{\lambda_n}(|\beta_nj|),\]
where the loss function is least absolute deviation(LAD for short), that does not need the noise obeys a gaussian distribution and be more
robust than least squares estimation. In fact, LAD estimation is the special case of M-estimation, which is named by Huber(1964, 1973, 1981)\cite{Huber1} \cite{Huber2} \cite{Huber3}firstly and can be obtained by minimizing the objective function
$$Q_n(\beta_n)=\frac{1}{n}\sum_{i=1}^{n}\rho(y_i-x_i^T\beta_n),$$
where the function $\rho$ can be selected.
For example, if we choose $\rho(x)=\frac{1}{2}x^21_{|x|\leq c}+(c|x|-c^2/2)1_{|x|>c}$,
where $c>0$, Huber estimator can be obtained;
if we choose $\rho(x)=|x|^q$, where $1\leq q\leq2$,
$L_q$ estimator will be obtained, with two special cases: LAD estimator for $q=1$ and OLS estimator for $q=2$.
If we choose
$\rho(x)=\alpha x^++(1-\alpha)(-x)^+$, where $0<\alpha<1,x^+=\max(x,0)$,
we call it quantile regression, and can also get LAD estimator for $\alpha=1/2$ especially.
\par When $p_n$ approaches infinity as $n$ tends to infinity, we assume that the function $\rho$ is convex and not monotone, and the monotone function $\varphi$ is the derivative of $\rho$. By imposing the appropriate regularity conditions,
Huber(1973), Portnoy(1984)\cite{Portnoy},Welsh(1989)\cite{Welsh} and Mammen(1989)\cite{Mammen} have proved that the M-estimator enjoyed the properties of consistency and asymptotic normality, where Welsh(1989) gave the weaker condition imposed on $\varphi$ and the stronger condition on $p_n/n$.
Bai and Wu \cite{Bai1} further pointed that the condition on $p_n$ could be a part of the integrable condition imposed on design matrix. Moreover, He and Shao(2000)\cite{He} studied the asymptotic properties of M-estimator in the case of the generalized model setting and the dimension $p_n$ getting bigger and bigger. Li(2011)\cite{Li} obtained the Oracle property of non-concave penalized M-estimator in high-dimensional model with the condition of $p_n\log n/n\rightarrow0,p_n^2/n\rightarrow0$, and proposed RSIS to make variable selection by applying rank sure independence screening method in the ultra high-dimensional model. Zou and Li(2008)\cite{Zou} combined penalized function and local linear approximation method(LLA) to prove that the obtained estimator enjoyed good asymptotic properties, and demonstrated this method improved the computational efficiency of local quadratic approximation(LQA) in the part of simulation.\\
Inspired by this, in this paper we consider the following problem:
\[\hat{\beta_n}=\text{arg}\min_{\beta_n} Q_n(\beta_n),\ Q_n(\beta_n)=\frac{1}{n}\sum_{i=1}^{n}\rho(y_i-x_i^T\beta_n)+\sum_{j=1}^{p_n}p'_{\lambda_n}(|\tilde{\beta}_{nj}|)|\beta_{nj}|,\tag{1.1}\]
where $p'_{\lambda_n}(\cdot)$ is the derivative of the penalized function, and $\tilde{\beta}_{n}=(\tilde{\beta}_{n1},\tilde{\beta}_{n2},...,\tilde{\beta}_{np_n})^T$ is the non-penalized estimator.
\par In this paper, we assume that the function $\rho$ is convex, hence the objective function is still convex and the obtained local minimizer is global minimizer.

\section*{\normalsize\bf 2 \quad{\normalsize\bf Main results}}
 For the convenience of statement, we first give some notations. Let
$\beta_0=(\beta_{01},\beta_{02},...,\beta_{0p})^T$ be the true parameter. Without loss of generality, we assume the first $k_n$ coefficients of covariates are nonzero, $p_n-k_n$ be coviariates with zero coefficients. $\beta_0=(\beta_{0(1)}^T,\beta_{0(2)}^T)^T,\hat{\beta_n}=(\hat{\beta}_{n(1)}^T,\hat{\beta}_{n(2)}^T)^T$ correspondingly. For the given symmetric matrix $Z$, denote by $\lambda_{min}(Z)$ and $\lambda_{max}(Z)$ the minimum and maximum eigenvalue of $Z$, respectively. Denote $\frac{X^TX}{n}:=D$ and
$D=\left(
\begin{array}{cc}
                                    D_{11} &D_{12}  \\
                                   D_{21}& D_{22} \\
                                  \end{array}
                                \right),$
where $D_{11}=\frac{1}{n}X_{(1)}^TX_{(1)}$. Finally we denote that $c_n=\max\{|p'_{\lambda_n}(|\tilde{\beta}_{nj}|)|:\tilde{\beta}_{nj}\neq0,1\leq j\leq p_n\}$.

\par Next, we state some assumptions which will be needed in the following results.\\
$(A_1)$  The function $\rho$ is convex on $R$, and its left derivative and right derivative $\varphi_+(\cdot),\ \varphi_-(\cdot)$ satisfies that
$\varphi_-(t)\leq \varphi(t)\leq\varphi_+(t),\ \forall t\in R$.\\
$(A_2)$ The error term $\varepsilon$ is i.i.d, and the distribution function $F$ of  $\varepsilon_i$ satisfies $F(S)=0$, where $S$ is the set of discontinuous points of $\varphi$.
 \par Moreover,
$E[\varphi(\varepsilon_i)]=0,\ 0<E[\varphi^2(\varepsilon_i)]=\sigma^2<\infty$, and
$G(t)\equiv E[\varphi(\varepsilon_i+t)]=\gamma t+o(|t|)$, where$\gamma>0$. Besides these, we assume that
$\displaystyle\lim_{t\rightarrow 0}E[\varphi(\varepsilon_i+t)-\varphi(\varepsilon_i)]^2=0$.\\
$(A_3)$ There exist constants $\tau_{1},\ \tau_{2},\ \tau_{3},\ \tau_{4}$ such that $0<\tau _{1}\leq \lambda_{min}(D)\leq\lambda_{max}(D)\leq\tau_{2}$ and
$0<\tau _{3}\leq \lambda_{min}(D_{11})\leq\lambda_{max}(D_{11})\leq\tau_{4}$.\\
$(A_4)$ $\lambda_{n}\rightarrow0(n\rightarrow\infty)$,\ 且$p_n=O(n^{1/2}),c_n=O(n^{-1/2})$.\\
$(A_5)$ Let $z_i$ be the transpose of the $i$th row vector of $X_{(1)}$, such that
$\displaystyle\lim_{n\rightarrow \infty}n^{-\frac{1}{2}}\max_{1\leq i\leq n}z_i^Tz_i=0.$
\par It is worth mentioning that conditions $(A_1)$ and $(A_2)$ are classical assumptions for M-estimation in linear model, which can be found in many references, for example Bai, Rao and Wu(1992)\cite{Bai}and Wu(2007)\cite{Wu}. The condition $(A_3)$ is frequently used for sparse model in the linear model regression theory, which requires that the eigenvalues of  the matrices $D$ and $D_{11}$ are bounded. The condition $(A_4)$ is weaker than that in previous references. In the condition $(A_4)$ we broad the order of $p_n$ to $n^{1/2}$, but in the references Huber(1973) and Li,Peng and Zhu(2011)\cite{Li} they required that $p_n^2/n\rightarrow0$, Portnoy(1984) required $p_n\log p_n/n\rightarrow0$,
and Mammen(1989) required $p_n^{3/2}\log p_n/n\rightarrow0$. Compared with these results, it is obvious that our sparse condition is much weaker. The condition $(A_5)$  is the same as that in Huang, Horowitz and Ma(2008)\cite{Huang}, which is used to prove the asymptotic properties of the nonzero part of M-estimation.

\textbf{Theorem 2.1}\ \textbf{(Consistency of estimator)} If the conditions $(A_1)-(A_4)$ hold, there exists a non-concave penalized M-estimation $\hat{\beta_n}$, such that \[\|\hat{\beta_n}-\beta_{0}\|=O_P((p_n/n)^{1/2}).\]

\textbf{Remark 2.1} From Theorem 2.1, we can obtain that there exists a global M-estimation $\hat{\beta_n}$ if we choose the appropriate tuning parameter $\lambda_{n}$, moreover this M-estimation is $(n/p_n)^{1/2}$-consistent. This convergence rate is the same as that in the references  Huber(1973) and Li,Peng and Zhu(2011).

\textbf{Theorem 2.2}\ \textbf{(The sparse of the model)} If the conditions $(A_1)-(A_4)$ hold and $\lambda_{min}(D)>\lambda_{max}(\frac{1}{n}\sum_{i=1}^{n}J_iJ_i^T)$, for the non-concave penalized M-estimation $\hat{\beta_n}$
we have \[P(\hat{\beta}_{n(2)}=0)\rightarrow 1.\]

\textbf{Remark 2.2} By Theorem 3.2, we can get that under the suitable conditions the global M-estimation of zero-coefficient variables goes to zero with a high probability when $n$ is large enough. This also shows that the model is sparse.\\

\textbf{Theorem 2.3}\ {\bf(Oracle property)} If the conditions $(A_1)-(A_5)$ hold and $\lambda_{min}(D)>\lambda_{max}(\frac{1}{n}\sum_{i=1}^{n}J_iJ_i^T)$,
with probability converging to one the non-concave penalized M-estimation
$\hat{\beta_n}=(\hat{\beta}_{n(1)}^T,\hat{\beta}_{n(2)}^T)^T$
has the following properties:
\par (1)(The consistency of the model selection)$\hat{\beta}_{n(2)}=0$;
\par (2)(Asymptotic normality)
\[\aligned\sqrt{n}s_n^{-1}u^T(\hat{\beta}_{n(1)}-\beta_{0(1)})&=\sum_{i=1}^{n}n^{-1/2}s_n^{-1}\gamma^{-1}u^TD_{11}z_i^T\varphi(\varepsilon_i)+o_P(1)\\
&\stackrel{\mbox{\small d}}{\longrightarrow}N(0,1),\endaligned\]
where $s_n^2=\sigma^2\gamma^{-1}u^TD_{11}^{-1}u$, and $u$ is any $k_n$ dimensional vector such that $\|u\|\leq1$. Meanwhile, $z_i$ is the transpose of the $i$th row vector of a $k_n\times k_n$ matrix $X_{(1)}$.

\textbf{Remark 2.3} From Theorem 2.3, M-estimation enjoys Oracle property, that is, the adaptive bridge estimator can
correctly select covariates with nonzero coefficients with probability converging to one
and that the estimator of nonzero coefficients has the same asymptotic distribution
that they would have if the zero coefficients were known in advance.

\textbf{Remark 2.4} In Fan and Peng(2004)\cite{Fan}, the authors obtained that the non-concave penalized M-estimation has the property of consistency with the condition $p_n^4/n\rightarrow0$,
and enjoyed the property of asymptotic normality with the condition $p_n^5/n\rightarrow0$. By Theorem 3.1-3.3, we can see that the corresponding conditions we exert is quite weak.\\

\section*{\normalsize\bf 3 \quad{\normalsize\bf Proofs of main results}}
\textbf{The proof of Theorem 2.1:}\
Let $\alpha_n=(p_n/n)^{1/2}+p_n^{1/2}c_n$,
where $u$ is a any $p_n$-dimensional vector such that $\|u\|=C$.
In the following part we only need to prove that
there exists a great enough positive constant $C$ such that
\[\liminf_{n\rightarrow\infty}P\{\inf_{\|u\|=C}Q_n(\beta_0+\alpha_n u)>Q_n(\beta_0)\}\geq 1-\varepsilon, \tag{3.1}\]
for any $\varepsilon>0$, that is,
there at least exists a local minimizer $\hat{\beta_n}$ such that $\|\hat{\beta_n}-\beta_{0}\|=O_P(\alpha_n)$
in the closed ball $\{\beta_0+\alpha_n u:\|u\|\leq C\}$.
Firstly by the triangle inequality we can get that
\begin{equation}
\aligned &Q_n(\beta_0+\theta u)-Q_n(\beta_0)\\
&=\frac{1}{n}\sum_{i=1}^{n}[\rho(y_i-x_i^T(\beta_0+\alpha_n u))-\rho(y_i-x_i^T\beta_0)]+\sum_{j=1}^{p_n}p'_{\lambda_n}(|\tilde{\beta}_{nj}|)(|\beta_{0j}+\alpha_n u_j|-|\beta_{0j}|)\\
&\geq \frac{1}{n}\sum_{i=1}^{n}[\rho(y_i-x_i^T(\beta_0+\alpha_n u))-\rho(y_i-x_i^T\beta_0)]-\alpha_n\sum_{j=1}^{p_n}p'_{\lambda_n}(|\tilde{\beta}_{nj}|)|u_j|\\
&:=T_1+T_2,
\endaligned\tag{3.2}
\end{equation}
where $T_1=\frac{1}{n}\sum_{i=1}^{n}[\rho(y_i-x_i^T(\beta_0+\alpha_n u))-\rho(y_i-x_i^T\beta_0)]$, $T_2=-\alpha_n\sum_{j=1}^{p_n}p'_{\lambda_n}(|\tilde{\beta}_{nj}|)|u_j|$.
Noticing that
\begin{equation}
\aligned T_1&=\frac{1}{n}\sum_{i=1}^{n}[\rho(y_i-x_i^T(\beta_0+\alpha_n u))-\rho(y_i-x_i^T\beta_0)]\\
&=\frac{1}{n}\sum_{i=1}^{n}[\rho(\varepsilon_i-\alpha_n x_i^Tu)-\rho(\varepsilon_i)]\\
&=\frac{1}{n}\sum_{i=1}^{n}\int_0^{-\alpha_n x_i^Tu}[\varphi(\varepsilon_i+t)-\varphi(\varepsilon_i)]dt-\frac{1}{n}\alpha_n \sum_{i=1}^{n}\varphi(\varepsilon_i)x_i^Tu\\
&:=T_{11}+T_{12},\endaligned\tag{3.3}
\end{equation}
where $T_{11}=\frac{1}{n}\sum_{i=1}^{n}\int_0^{-\alpha_n x_i^Tu}[\varphi(\varepsilon_i+t)-\varphi(\varepsilon_i)]dt$, $T_{12}=-\frac{1}{n}\alpha_n \sum_{i=1}^{n}\varphi(\varepsilon_i)x_i^Tu$,\\
combining with Von-Bahr Esseen inequality and the fact that $|T_{12}|\leq\frac{1}{n}\alpha_n \|u\|\|\sum_{i=1}^{n}\varphi(\varepsilon_i)x_i\|$, we instantly have
\[E[\|\sum_{i=1}^{n}\varphi(\varepsilon_i)x_i\|^2]
\leq n\sum_{i=1}^{n}E[\|\varphi(\varepsilon_i)x_i\|^2]=n\sum_{i=1}^{n}E[\varphi^2(\varepsilon_i)x_i^Tx_i
\leq n^2p_n\sigma^2,\]
hence
\[|T_{12}|=O_P(\alpha_np_n^{1/2})\|u\|=O_P((p_n^2/n)^{1/2}).\tag{3.4}\]
Secondly for $T_{11}$, let $T_{11}=\sum_{i=1}^{n}A_{in}$,
where $A_{in}=\frac{1}{n}\int_0^{-\alpha_n x_i^Tu}[\varphi(\varepsilon_i+t)-\varphi(\varepsilon_i)]dt$,
so
$$T_{11}=\sum_{i=1}^{n}[A_{in}-E(A_{in})]+\sum_{i=1}^{n}E(A_{in}):=T_{111}+T_{112}.$$
We can easily obtain
$E(T_{111})=0$.
From Von-Bahr Esseen inequality, Schwarz inequality and the condition $(B_3)$,
it follows that
\[\aligned var(T_{111})&=var(\sum_{i=1}^{n}A_{in})\leq \frac{1}{n}\sum_{i=1}^{n}E\biggl(\int_0^{-\alpha_n x_i^Tu}[\varphi(\varepsilon_i+t)-\varphi(\varepsilon_i)]dt\biggr)^2\\
&\leq \frac{1}{n}\sum_{i=1}^{n}|\alpha_n x_i^Tu||\int_0^{-\alpha_n x_i^Tu}E[\varphi(\varepsilon_i+t)-\varphi(\varepsilon_i)]^2dt|\\
&=\frac{1}{n}\sum_{i=1}^{n}o_P(1)(\alpha_n x_i^Tu)^2=\frac{1}{n}o_P(1)\alpha_n^2\sum_{i=1}^{n}u^Tx_ix_i^Tu\\
&=o_P(1)\alpha_n^2u^TDu\leq \lambda_{max}(D)o_P(1)\alpha_n^2\|u\|^2=o_P(\alpha_n^2)\|u\|^2,
\endaligned\]
together by Markov inequality yields that
$$P(|T_{111}|>C_1\alpha_n\|u\|)\leq \frac{var(T_{111})}{C_1^2\alpha_n^2\|u\|^2}\leq\frac{o_P(\alpha_n^2)\|u\|^2}{C_1^2\alpha_n^2\|u\|^2}\rightarrow0(n\rightarrow\infty),$$
hence \[T_{111}=o_P(\alpha_n)\|u\|.\tag{3.5}\]
As for $T_{112}$,
\[\aligned T_{112}&=\sum_{i=1}^{n}E(A_{in})=\frac{1}{n}\sum_{i=1}^{n}\int_0^{-\alpha_n x_i^Tu}[\gamma t+o(|t|)]dt\\
&=\frac{1}{n}\sum_{i=1}^{n}(\frac{1}{2}\gamma\alpha_n^2u^Tx_ix_i^Tu+o_P(1)\alpha_n^2u^Tx_ix_i^Tu)\\
&=\frac{1}{2}\gamma\alpha_n^2u^TDu+o_p(1)\alpha_n^2u^TDu\\
&\geq[\frac{1}{2}\gamma\lambda_{min}(D)+o_P(1)]\alpha_n^2\|u\|^2.
\endaligned\tag{3.6}\]
Finally considering $T_2$, we can easily obtain
\[T_{2}\leq (p_n)^{1/2}\alpha_n\max\{|p'_{\lambda_n}(|\tilde{\beta}_{nj}|)|,1\leq j\leq k_n\}\|u\|
=(p_n)^{1/2}\alpha_nc_n\|u\|\leq \alpha_n^2\|u\|.
\tag{3.7}\]
This together with (3.3)-(3.7) yields that we can choose a great enough constant $C$  such that $T_{111}$ and $T_2$ is controlled by $T_{112}$, which follows that there at least exists a local minimizer $\hat{\beta_n}$ such that $\|\hat{\beta_n}-\beta_{0}\|=O_P(\alpha_n)$
in the closed ball $\{\beta_0+\alpha_n u:\|u\|\leq C\}$. $\Box$\\

\textbf{The proof of Theorem 2.2:}\ From Theorem 2.1, as long as we choose a great enough constant $C$ and appropriate $\alpha_n$, then $\hat{\beta_n}$ will be in the ball $\{\beta_0+\alpha_n u:\|u\|\leq C\}$ with probability converging to one, where
\ $\alpha_n=(p_n/n)^{1/2}+p_n^{1/2}c_n$.
For any $p_n$-dimensional vector $\beta_n$, now we denote
$\beta_n=(\beta_{n(1)}^T,\beta_{n(2)}^T)^T,\ \beta_{n(1)}=\beta_{0(1)}+\alpha_n u_{(1)},\ \beta_{n(2)}=\beta_{0(2)}+\alpha_n u_{(2)}=\alpha_n u_{(2)}$, where $\beta_0=(\beta_{0(1)}^T,\beta_{0(2)}^T)^T,\ \|u\|^2=\|u_{(1)}\|^2+\|u_{(2)}\|^2\leq C^2$.
Meanwhile let $$V_n(u_{(1)},u_{(2)})=Q_n((\beta_{n(1)}^T,\beta_{n(2)}^T)^T)-Q_n((\beta_{0(1)}^T,0^T)^T),$$
then by minimizing $V_n(u_{(1)},u_{(2)})$ we can obtain the estimator $\hat{\beta_n}=(\hat{\beta}_{n(1)}^T,\hat{\beta}_{n(2)}^T)^T$, where $\|u_{(1)}\|^2+\|u_{(2)}\|^2\leq C^2$.
In the following part, we will prove that as long as $\|u\|\leq C,\ \|u_{(2)}\|>0$,
\[P(V_n(u_{(1)},u_{(2)})-V_n(u_{(1)},0)>0)\rightarrow1(n\rightarrow\infty)\tag{3.8}\]
holds, for any $p_n$-dimensional vector $u=(u_{(1)}^T,u_{(2)}^T)^T$.
We can easily find the fact that
\[\aligned &V_n(u_{(1)},u_{(2)})-V_n(u_{(1)},0)
=Q_n((\beta_{n(1)}^T,\beta_{n(2)}^T)^T)-Q_n((\beta_{n(1)}^T,0^T)^T)\\
&=\frac{1}{n}\sum_{i=1}^{n}[\rho(\varepsilon_i-\alpha_n H_i^Tu_{(1)}-\alpha_n J_i^Tu_{(2)})-\rho(\varepsilon_i-\alpha_n H_i^Tu_{(1)})]
+\sum_{j=k_n+1}^{p_n}p'_{\lambda_n}(|\tilde{\beta}_{nj}|)|\alpha_n u_j|\\
&=\frac{1}{n}\sum_{i=1}^{n}\int_{-\alpha_n H_i^Tu_{(1)}}^{-\alpha_n H_i^Tu_{(1)}-\alpha_n J_i^Tu_{(2)}}
[\varphi(\varepsilon_i+t)-\varphi(\varepsilon_i)]dt-\frac{1}{n}\alpha_n\sum_{i=1}^{n}\varphi(\varepsilon_i)J_i^Tu_{(2)}\\
&\ \ +\sum_{j=k_n+1}^{p_n}p'_{\lambda_n}(|\tilde{\beta}_{nj}|)|\alpha_n u_j|:=W_1+W_2+W_3,
\endaligned\tag{3.9}\]
where $H_i$ and $J_i$ are $k_n$ and $p_n-k_n$ dimensional vectors respectively such that $x_i=(H_i^T+J_i^T)^T$.
Similar to the proof of Theorem 3.1, we get that
\[\aligned W_1&=\frac{1}{n}\sum_{i=1}^{n}\int_{-\alpha_n H_i^Tu_{(1)}}^{-\alpha_n H_i^Tu_{(1)}-\alpha_n J_i^Tu_{(2)}}
[\varphi(\varepsilon_i+t)-\varphi(\varepsilon_i)]dt\\
&=\frac{1}{2n}\sum_{i=1}^{n}\gamma\alpha_n^2u^Tx_ix_i^Tu-\frac{1}{2n}\sum_{i=1}^{n}\gamma\alpha_n^2u_{(2)}^TJ_iJ_i^Tu_{(2)}
+o_P(1)\alpha_n^2\|u\|^2+o_P(1)\alpha_n\|u\|\\
&\geq\frac{1}{2}\gamma\alpha_n^2[\lambda_{min}(D)-\lambda_{max}(\frac{1}{n}\sum_{i=1}^{n}J_iJ_i^T)]\|u\|^2
+o_P(1)\alpha_n^2\|u\|^2+o_P(1)\alpha_n\|u\|，
\endaligned\tag{3.10}\]
\[|W_2|=|-\frac{1}{n}\alpha_n\sum_{i=1}^{n}\varphi(\varepsilon_i)J_i^Tu_{(2)}|
=O_P((p_n^2/n)^{1/2})\|u\|,\tag{3.11}\]
and
\[\aligned|W_3|&=|\sum_{j=k_n+1}^{p_n}p'_{\lambda_n}(|\tilde{\beta}_{nj}|)|\alpha_n u_j||
\leq (p_n)^{1/2}\alpha_n\max\{|p'_{\lambda_n}(|\tilde{\beta}_{nj}|)|,k_n+1\leq j\leq p_n\}\|u\|\\
&=(p_n)^{1/2}\alpha_nc_n\|u\|\leq \alpha_n^2\|u\|.\endaligned\tag{3.12}\]
By formula (3.10)-(3.12) and the condition $\lambda_{min}(D)>\lambda_{max}(\frac{1}{n}\sum_{i=1}^{n}J_iJ_i^T)$, it follows that
\[\aligned &V_n(u_{(1)},u_{(2)})-V_n(u_{(1)},0)
\geq\frac{1}{2}\gamma\alpha_n^2[\lambda_{min}(D)-\lambda_{max}(\frac{1}{n}\sum_{i=1}^{n}J_iJ_i^T)]\|u\|^2\\
& +o_P(1)\alpha_n^2\|u\|^2+o_P(1)\alpha_n\|u\|]+O_P((p_n^2/n)^{1/2})\|u\|+O_P(\alpha_n^2)\|u\|>0,
\endaligned\]
which yields that as long as $\|u\|\leq C,\ \|u_{(2)}\|>0$,
\[P(V_n(u_{(1)},u_{(2)})-V_n(u_{(1)},0)>0)\rightarrow1(n\rightarrow\infty)\]
holds, for any $p_n$-dimensional vector $u=(u_{(1)}^T,u_{(2)}^T)^T$. $\Box$

\textbf{The proof of Theorem 2.3:}\ \ It is obvious that the conclusion (1) can be obtained instantly by Theorem 2.2, so we only need to prove the conclusion (2).
It follows from Theorem 2.1 that $\hat{\beta_n}$ is consistent of $\beta_0$
and $\hat{\beta}_{n(2)}=0$ with probability converging to one from Theorem 2.2.
Therefore $\hat{\beta_{n(1)}}$ holds that
\[\frac{\partial Q_n(\beta_n)}{\partial\beta_{n(1)}}\mid_{\beta_{n(1)}=\hat{\beta}_{n(1)}}=0,\]
that is
\[-\frac{1}{n}\sum_{i=1}^{n}H_i\varphi(y_i-H_i^T\hat{\beta}_{n(1)})+W_{(1)}=0,\tag{3.13}\]
where \[W=(p'_{\lambda_n}(|\tilde{\beta}_{n1}|)sgn(\hat{\beta}_{n1}),p'_{\lambda_n}(|\tilde{\beta}_{n2}|)sgn(\hat{\beta}_{n2}),
\cdots,p'_{\lambda_n}(|\tilde{\beta}_{np_n}|)sgn(\hat{\beta}_{np_n}))^T.\]
In the following part we give the Taylor expansion of upper left first term:
\[\aligned -\frac{1}{n}\sum_{i=1}^{n}\{H_i\varphi(y_i-H_i^T\hat{\beta}_{0(1)})-
[\varphi'(y_i-H_i^T\beta_{0(1)})H_iH_i^T+o_P(1)](\hat{\beta}_{n(1)}-\beta_{0(1)})\}+W_{(1)}=0.\endaligned\]
Noticing that $y_i=H_i^T\beta_{0(1)}+\varepsilon_i$, we have
\[\aligned -\frac{1}{n}\sum_{i=1}^{n}H_i\varphi(\varepsilon_i)+\frac{1}{n}\sum_{i=1}^{n}[\varphi'(\varepsilon_i)H_iH_i^T+o_P(1)]
(\hat{\beta}_{n(1)}-\beta_{0(1)})+W_{(1)}=0,\endaligned\]
which yields that
\[\aligned \frac{1}{n}\gamma\sum_{i=1}^{n}H_iH_i^T(\hat{\beta}_{n(1)}-\beta_{0(1)})
&=\frac{1}{n}\sum_{i=1}^{n}H_i\varphi(\varepsilon_i)-W_{(1)}+(\hat{\beta}_{n(1)}-\beta_{0(1)})o_P(1)\\
&+\frac{1}{n}\sum_{i=1}^{n}(\gamma-\varphi'(\varepsilon_i))H_iH_i^T(\hat{\beta}_{n(1)}-\beta_{0(1)}).
\endaligned\]
Then as long as $\|u\|\leq 1$,
\[\aligned u^T(\hat{\beta}_{n(1)}-\beta_{0(1)})
&=n^{-1}\gamma^{-1}u^TD_{11}^{-1}\sum_{i=1}^{n}H_i\varphi(\varepsilon_i)\\
&+n^{-1}\gamma^{-1}u^TD_{11}^{-1}\sum_{i=1}^{n}(\gamma-\varphi'(\varepsilon_i))H_iH_i^T(\hat{\beta}_{n(1)}-\beta_{0(1)})\\
&-\gamma^{-1}u^TD_{11}^{-1}W_{(1)}+o_P(\alpha_n)
\endaligned\tag{3.14}\]
holds, for any $k_n$-dimensional vector $u$.
For upper right third term, we can obtain
\[\aligned |\gamma^{-1}u^TD_{11}^{-1}W_{(1)}|\leq\frac{1}{\gamma\lambda_{min}(D_{11})}\|W_{(1)}\|
\leq\frac{1}{\gamma\lambda_{min}(D_{11})}p_n^{1/2}c_n
\leq\frac{\alpha_n}{\gamma\lambda_{min}(D_{11})}\rightarrow o_P(1)(n\rightarrow\infty)
.\endaligned\tag{3.15}\]
Now let us deal with upper right second term. Theorem 2.1 and the condition $(A_3)$ yield that
\[\aligned &|n^{-1}\gamma^{-1}u^TD_{11}^{-1}\sum_{i=1}^{n}(\gamma-\varphi'(\varepsilon_i))H_iH_i^T(\hat{\beta}_{n(1)}-\beta_{0(1)})|\\
&\leq\frac{1}{n\gamma\lambda_{min}(D_{11})}\|\sum_{i=1}^{n}(\gamma-\varphi'(\varepsilon_i))H_iH_i^T(\hat{\beta}_{n(1)}-\beta_{0(1)})\|\\
&\leq\frac{1}{n\gamma\lambda_{min}(D_{11})}\|\sum_{i=1}^{n}(\gamma-\varphi'(\varepsilon_i))H_iH_i^T\|\|\hat{\beta}_{n(1)}-\beta_{0(1)}\|\\
&\leq\frac{O_P(1)}{n\gamma\lambda_{min}(D_{11})}\|\hat{\beta}_{n(1)}-\beta_{0(1)}\|=O_P(p_n^{1/2}n^{-3/2}),
\endaligned\tag{3.16}\]
where the upper third inequality sign holds because of applying Lemma 3 of Mammen(1989).
Combining (3.15)-(3.17), we have
\[\aligned u^T(\hat{\beta}_{n(1)}-\beta_{0(1)})
=n^{-1}\gamma^{-1}u^TD_{11}^{-1}\sum_{i=1}^{n}H_i\varphi(\varepsilon_i)+O_P(\alpha_n)+O_P(p_n^{1/2}n^{-3/2})，
\endaligned\]
that is,
\[n^{1/2}u^T(\hat{\beta}_{n(1)}-\beta_{0(1)})=n^{-1/2}\gamma^{-1}u^TD_{11}^{-1}\sum_{i=1}^{n}H_i\varphi(\varepsilon_i)+o_P(1).\tag{3.17}\]
Denote $s_n^2=\sigma^2\gamma^{-1}u^TD_{11}^{-1}u,\ F_{in}=n^{-1/2}s_n^{-1}\gamma^{-1}u^TD_{11}^{-1}z_i^T$,
where $z_i$ is a $k_n\times k_n$ matrix and  the transpose of the $i$th row vector of $X_{(1)}$,
then
$n^{1/2}u^T(\hat{\beta}_{n(1)}-\beta_{0(1)})=\sum_{i=1}^{n}F_{in}\varphi(\varepsilon_i)+o_P(1)$.
It follows from $(B_5)$ that
\[\aligned \sum_{i=1}^{n}F^2_{in}&=\sum_{i=1}^{n}F_{in}F'_{in}
=\sum_{i=1}^{n}(n^{-1/2}s_n^{-1}\gamma^{-1}u^TD_{11}^{-1}z_i^T)(n^{-1/2}s_n^{-1}\gamma^{-1}z_iD_{11}^{-1}u)\\
&=\sum_{i=1}^{n}n^{-1}s_n^{-2}\gamma^{-2}u^TD_{11}^{-1}z_i^Tz_iD_{11}^{-1}u
=s_n^{-2}\gamma^{-2}u^TD_{11}^{-1}u=\sigma^{-2}.
\endaligned\]
Applying Slutsky Theorem, we obtain that
\[\sqrt{n}s_n^{-1}u^T(\hat{\beta}_{n(1)}-\beta_{0(1)})\stackrel{\mbox{\small d}}{\longrightarrow}N(0,1).\;\Box\]

\section*{\normalsize\bf 4 \quad{\normalsize\bf Simulation results} }
\par In this section we evaluate the performance of the M-estimator proposed in (1.1) by simulation studies.
\par About the data. Simulate data by the model $Y=X\beta+\varepsilon$, where $\beta_{0(1)}=(-2,2.5,3,-1)^T$,
$\varepsilon$ follows $N(0,1), t_5$
and mixed normally distribution $0.9N(0,1)+0.1N(0,9)$ respectively. And the design matrix $X$ is generated by p-dimensional multivariate normal distribution with mean zero and covariance matrix whose $(i,j)$th component component is $\rho^{|i-j|}$, where we set $\rho=0.5$.
\par About loss function. In this section we can choose some special loss functions, such as LAD loss function, OLS loss function and Huber loss function. In this paper we choose LAD loss function.
\par About penalty function. For $p'_{\lambda_n}(|\widetilde{\beta}_{nj}|)$ in penalty function, we choose penalty function as SACD estimation in the following:
$$p_{\lambda_n}(|\beta|)=\left\{%
                           \begin{array}{ll}
                             \lambda_n|\beta|, & 0\leq|\beta|\leq\lambda_n, \\
                             -(\beta^2-2a\lambda_n|\beta|+\lambda_n^2)/(2(a-1)), & \lambda_n<|\beta|< a\lambda_n, \\
                            (a+1)\lambda_n^2/2, & |\beta|> a\lambda_n,\\
                           \end{array}%
                         \right.
$$
then
$p'_{\lambda_n}(|\widetilde{\beta}_{nj}|)=
\lambda_nI(|\widetilde{\beta}_{nj}|\leq\lambda_n)+\frac{a\lambda_n-|\widetilde{\beta}_{nj}|}{a-1}I(\lambda_n<|\widetilde{\beta}_{nj}|\leq a\lambda_n)$. By the proposal of Fan and Li(2001), we can select $a=3.7$, which yields that
generalized cross validation can be applied in searching the best tuning parameter $\lambda_n$.
\par About stimulation algorithm. For proposed LLA method, we connect penalty function with independent variables and independent variable respectively, then programme by using quantile package in R. For Lasso method, we use Lars package to simulate.
\par About the selection of tuning parameter. We apply BIC criterion to select tuning parameter. The criterion is in the following
$$
BIC(\lambda_n)=\ln(\frac{1}{n}\sum_{i=1}^{n}\rho(y_i-x_i^T\hat{\beta}))+DF_{\lambda_n}\ln(n)/n,
$$
where $DF_{\lambda_n}$ is the generalized degree of freedom in the reference Fan and Li(2001).
\par About selection of evaluation index. In order to evaluate the performance of the  estimators, we select four measures called EE, PE, C, IC and CP which are obtained by 500 replicates. EE is median of $||\hat{\beta}-\beta_0||_2$ to evaluate the estimation accuracy, and PE is the prediction error defined by median of $n^{-1}\|Y-X\hat{\beta}\|^2$. The other three measures are to qualify the performance of model consistency, where C and IC refer to the average number of correctly selected zero covariates and the average number of incorrectly selected zero covariates, and CP is the proportion of the number of the correct selection of zero variables to the total number of zero variables.
\par In the following we will compare the performances of the method LLA we proposed, Lasso method and Oracle estimation. Set $n=200,500,700$ respectively and $p=[2\sqrt{n}]$.

\par From table 1, we notice that the index EE, C, IC, CP of our proposed method LLA perform better when $\varepsilon\sim N(0,1)$. In particular, for the index CP, LLA outperforms Lasso. The reason of this may be that we impose different penalties for important and unimportant variables, while Lasso imposes the same penalties for all variables. Moreover, with the increase of sample size, the ability of LLA method to correctly identify unimportant variables is also increasing. When the sample size is 700 and the number of explanatory variables is 53, an average of 48.9617 unimportant variables-zero variables are estimated to be zero on average, with an average accuracy of 99.92\%.
\par An interesting fact can be found from Table 2, that is, when the error term is chosen as $t_5$, the accuracy of the method LLA proposed to correctly exclude incorrect variables is slightly higher than that of the case where the error term is standardized normal distribution. The reason is that when the error term is heavy-tailed, it is more appropriate to choose LLA, but the accuracy of estimation and prediction is slightly worse than that of Lasso. When the sample size increases, the LLA and Oracle estimates perform equally well in the selection of important variables and the complexity of the model.
\par As can be seen from Table 3, when the error term is set to a mixed normal distribution, the ability of the proposed method to correctly select zero variables is good. In the case of small sample size, the ability of the Lasso method to select important variables is better.
%表格
\begin{table}[tbp]\footnotesize
\centering  % 表居中
\caption{ Simulation results for $\varepsilon\sim N(0,1)$ .}
\begin{tabular}{lcccccc}  % {lccc} 表示各列元素对齐方式，left-l,right-r,center-c
\hline
Setting &Method &EE &PE &C &IC &CP \\ \hline  % \hline 在此行下面画一横线
n=200&Oracle &10.8544 &3.3916 &24.0000 &0 &100\%\\
% \\表示重新开始一行
p=28&Lasso &10.5726	&3.3035	&10.8480 &0 &45.20\%\\
m=24&LLA &10.9153	&3.3947	&23.8540	&0	&99.39\%\\[2mm]
n=500 &Oracle	&19.9085	&5.4118	&41.0000	&0	&100\%\\
% \\ 表示重新开始一行
p=45&Lasso	&19.5952	&5.2928	&18.9920	&0	&46.32\%\\
m=41&LLA	&19.9233	&5.4045	&40.9140	&0	&99.79\%\\[2mm]
n=700&Oracle	&24.3006	&6.3847	&49.0000	&0	&100\%\\
  % \\ 表示重新开始一行
p=53&Lasso	&24.0315	&6.2994	&23.1009	&0	&47.14\%\\
m=49&LLA	&24.3666	&6.4077	&48.9617	&0	&99.92\%\\ \hline
\end{tabular}
\end{table}
%结束

%表格
\begin{table}[tbp]\footnotesize
\centering  % 表居中
\caption{ Simulation results for $\varepsilon\sim t_5$ .}
\begin{tabular}{lcccccc}  % {lccc} 表示各列元素对齐方式，left-l,right-r,center-c
\hline
Setting &Method &EE &PE &C &IC &CP \\ \hline  % \hline 在此行下面画一横线
n=200&Oracle &10.5634	&4.2892	&24.0000	&0	&100\%\\
p=28&Lasso &10.2810	&4.1649	&11.7700	&0	&49.04\%\\
m=24&LLA &10.6448	&4.2725	&23.8780	&0	&99.49\%\\[2mm]
n=500 &Oracle	&19.4296	&6.8240	&41.0000	&0	&100\%\\
p=45&Lasso	&19.1157	&6.7042	&18.9580	&0	&46.24\%\\
m=41&LLA	&19.4665	&6.8335	&40.9560	&0	&99.89\%\\[2mm]
n=700&Oracle	&23.7784	&8.0637	&49.0000	&0	&100\%\\
p=53&Lasso	&23.4389	&7.9551	&22.8800	&0	&46.69\%\\
m=49&LLA	&23.7808	&8.0919	&48.9740	&0	&99.94\%\\ \hline
\end{tabular}
\end{table}
%结束

%表格
\begin{table}[tbp]\footnotesize
\centering  % 表居中
\caption{ Simulation results for $\varepsilon\sim 0.9N(0,1)+0.1N(0,9)$ .}
\begin{tabular}{lcccccc}  % {lccc} 表示各列元素对齐方式，left-l,right-r,center-c
\hline
Setting &Method &EE &PE &C &IC &CP \\ \hline
n=200&Oracle &10.4815	&4.4830	&24.0000	&0	&100\%\\
p=28&Lasso &10.2030	&4.4063	&11.6360	&0	&48.48\%\\
m=24&LLA &10.5826	&4.4529	&23.9240	&0	&99.68\%\\[2mm]
n=500 &Oracle	&19.2539	&7.1997	&41.0000	&0	&100\%\\
p=45&Lasso	&18.9670	&7.0960	&19.3840	&0	&47.28\%\\
m=41&LLA	&19.2950	&7.1173	&40.9520	&0	&99.88\%\\[2mm]
n=700&Oracle	&23.6354	&8.5657	&49.0000	&0	&100\%\\
p=53&Lasso	&23.2424	&8.4609	&23.0580	&0	&47.06\%\\
m=49&LLA	&23.6566	&8.3699	&48.9300	&0	&99.86\%\\ \hline
\end{tabular}
\end{table}
%结束

%--------------------------the references---------------------------------------------------
\newpage
%\makeatletter
%\renewcommand\@biblabel[1]{}
%\renewenvironment{thebibliography}[1]
%     {\section*{\refname}%
%      \@mkboth{\MakeUppercase\refname}{\MakeUppercase\refname}%
%      \list{\@biblabel{\@arabic\c@enumiv}}%
%           {\settowidth\labelwidth{\@biblabel{#1}}%
%            \leftmargin\labelwidth
%            \advance\leftmargin\labelsep
%            \itemindent-1.5em
%            \@openbib@code
%            \usecounter{enumiv}%
%            \let\p@enumiv\@empty
%            \renewcommand\theenumiv{\@arabic\c@enumiv}}%
%      \sloppy
%      \clubpenalty4000
%      \@clubpenalty \clubpenalty
%      \widowpenalty4000%
%      \sfcode`\.\@m}
%     {\def\@noitemerr
%       {\@latex@warning{Empty `thebibliography' environment}}%
%      \endlist}
%\makeatother

%--------------------------the references---------------------------------------------------

%\end{CJK*}

{\footnotesize
\begin{thebibliography}{99}\baselineskip=6.mm
\setlength{\itemsep}{0pt} \setlength{\parskip}{0pt}
\bibitem{Huber1} Huber, P. {Robust estimation of a location parameter}, The Annals of Mathematical Statistics, \textbf{35}, 73-101，1964.
\bibitem{Huber2} Huber, P., {Robust regression: Asymptotics，conjectures and Monte Carlo}, The Annals of Statistics, \textbf{1}, 799-821, 1973.
\bibitem{Huber3} Huber, P., {Robust Statistics}, Wiley, 1981.
\bibitem{Portnoy} Portnoy, S., {Asymptotic behavior of M-estimators of p regression parameters when  $p^2/n$ is large，I：Consistency}, The Annals of Statistics, \textbf{12}, 1298-1309, 1984.
\bibitem{Welsh} Welsh, A., {On M-processes and M-estimation"，The Annals of Statistics}, \textbf{17}, 337-361, 1989.
\bibitem{Mammen}Mammen, E., {Asymptotics with increasing dimension for robust regression with applicationsto the bootstrap}, The Annals of Statistics, \textbf{17}, 382-400, 1989.
\bibitem{Bai1} Bai, Z., Wu, Y., {Limiting behavior of M-estimators of regression coefficients in high dimensional linear models I. scale-dependent case}, Journal of Multivariate Analysis, \textbf{51}, 211-239, 1994.
\bibitem{He}He, X., Shao, Q., {On parameters of increasing dimensions}, Journal of Multivariate Analysis, \textbf{73}, 120-135, 2000.
\bibitem{Li}Li, G., Peng, H., Zhu，L., {Nonconcave penalized M-estimation with a diverging number of parameters}, Statistica Sinica, \textbf{21}, 391-419, 2011.
\bibitem{Zou}Zou, H., Li, R., {One-step sparse estimates in nonconcave penalized likelihood models}, The Annals of Statistics, \textbf{36}, 1509-1566, 2008.
\bibitem{Bai} Bai, Z., Rao, C., Wu, Y., {M-estimation of multivariate linear regression parameters under a convex discrepancy function}, Statistica Sinica, \textbf{2}, 237-254, 1992.
\bibitem{Wu} Wu, W., {M-estimation of linear models with dependent errors}, The Annals of Statistics, \textbf{35}, 495-521, 2007.
\bibitem{Huang}Huang, J., Horowitz, J., Ma, S., {Asymptotic properties of bridge estimators in sparse high-dimensional regression models}, The Annals of Statistics, \textbf{36}, 587-613, 2008.
\bibitem{Fan} Fan, J., Peng, H., {Nonconcave penalized likelihood with a diverging number of parameters}, The Annals of Statistics, \textbf{32}, 928-961, 2004.
\end{thebibliography}
}
\end{document}